\documentclass[11pt]{article}
\usepackage[a4paper]{geometry}

\usepackage{amsmath}
\usepackage{amssymb}
\usepackage{amsthm}
\usepackage{mathrsfs}
\usepackage{bbm}
\usepackage{empheq}
\usepackage[loose]{subfigure}
\usepackage{epsfig}
\usepackage{graphicx}
\usepackage{psfrag}
\usepackage[usenames,dvipsnames]{pstricks}
\usepackage{pst-plot}
\usepackage[colorlinks]{hyperref}
\usepackage{tabls}
\usepackage{paralist}
\usepackage{hyperref}

\DeclareSymbolFontAlphabet{\Bbb}{AMSb}

\newlength{\fixboxwidth}
\newcommand{\COMMENT}[1]{}
\newcommand{\E}{\mathbb{E}}

\newcommand{\one}{\mathbbm{1}}
\renewcommand{\P}{\mathbb{P}}

\newcommand{\R}{\mathbb{R}}

\newcommand{\quark}{\setbox0\hbox{$x$}\hbox to\wd0{\hss$\cdot$\hss}}
\newcommand{\smid}{\,\middle|\,}

\newtheorem{thm}{Theorem}

\theoremstyle{definition}
\newtheorem{defn}{Definition}

\newtheorem{eg}{Example}
\newtheorem{pb}{Problem}

\hypersetup{
    linkcolor=blue,
}

\title{On the Brittleness of Bayesian Inference}

\author{Houman Owhadi$^{1,\ast}$, Clint Scovel$^{2}$, Tim Sullivan$^{3}$}

\date{\today}

\makeatletter
\@addtoreset{equation}{section}
\@addtoreset{figure}{section}
\@addtoreset{table}{section}
\makeatother

\renewcommand{\thefigure}{\arabic{figure}}

\makeatletter
\renewcommand{\p@subfigure}{\thefigure}
\makeatother

\newcounter{mycount}

\makeatletter
\def\blfootnote{\gdef\@thefnmark{}\@footnotetext}
\makeatother

\begin{document}

\maketitle

\begin{abstract}
	With the advent of high-performance computing, Bayesian methods are increasingly popular tools for the quantification of uncertainty throughout science and industry.
	Since these methods impact the making of sometimes critical decisions in increasingly complicated contexts, the sensitivity of their posterior conclusions with respect to the underlying models and prior beliefs is a pressing question for which there currently exist
	positive and negative results.
	We report new results suggesting that, although Bayesian methods are robust when the number of possible outcomes is finite or when only a finite number of marginals of the data-generating distribution are unknown, they could be generically brittle when applied to continuous systems (and their discretizations) with finite information on the data-generating distribution.
	If closeness is defined in terms of the total variation metric or the matching of a finite system of generalized moments, then (1) two practitioners who use arbitrarily close models and observe the same (possibly arbitrarily large amount of) data may reach opposite conclusions;  and (2) any given prior and model can be slightly perturbed to achieve any desired posterior conclusions.
	The mechanism causing brittlenss/robustness suggests that learning and robustness are antagonistic requirements and raises the question of a  missing stability condition for using Bayesian Inference in a continuous world under finite information.
\end{abstract}

\blfootnote{\noindent ${^1}$Department of  Computing + Mathematical Sciences, California Institute of Technology, Pasadena CA 91125, USA. E-mail: owhadi@caltech.edu}
\blfootnote{\noindent ${^2}$Department of  Computing + Mathematical Sciences, California Institute of Technology, Pasadena CA 91125, USA. E-mail: clintscovel@gmail.com}
\blfootnote{\noindent ${^3}$Mathematics Institute, University of Warwick, CV4 7AL, UK. E-mail: tim.sullivan@warwick.ac.uk}

\blfootnote{\noindent ${^\ast}$Corresponding author.} \setcounter{footnote}{3}

The application of Bayes' theorem in the form of Bayesian inference has fueled an ongoing debate with practical consequences in science, industry, medicine and law \cite{Efron:2013}.
One commonly-cited justification for the application of Bayesian reasoning is Cox's theorem \cite{Cox:1946}, which has been interpreted as stating that any `natural' extension of Aristotelian logic to uncertain contexts must be Bayesian \cite{Jaynes:2003}.
It has now been shown that Cox's theorem as originally formulated is incomplete \cite{Halpern:1999a} and there is debate about the `naturality' of the additional assumptions required for its validity \cite{Arnborg:2001, Dupre:2009, Halpern:1999b, Hardy:2002}, e.g.\ the assumption that knowledge can be always represented in the form of a $\sigma$-additive probability measure that assigns to each measurable event a \emph{single} real-valued probability.

However --- and this is the topic of this article --- regardless of the internal logic, elegance and appealing simplicity of Bayesian reasoning, a critical question is that of the robustness of its posterior conclusions with respect to perturbations of the underlying models and priors.

For example, a frequentist statistician might ask, if the data happen to be a sequence of i.i.d.\ draws from a fixed data-generating distribution $\mu^{\dagger}$, whether or not the Bayesian posterior will asymptotically assign full mass to a parameter value that corresponds to $\mu^{\dagger}$.
When it holds, this property is known as \emph{frequentist consistency} of the Bayes procedure, or the \emph{Bernstein--von Mises property}.

Alternatively, without resorting to a frequentist data-generating distribution $\mu^{\dagger}$, a Bayesian statistician who was also a numerical analyst might ask questions about stability and conditioning:  does the posterior distribution (or the posterior value of a particular quantity of interest) change only slightly when elements of the problem setup (namely the prior distribution, the likelihood model, and the observed data) are perturbed, e.g.\ as a result of observational error, numerical discretization, or algorithmic implementation?
When it holds, this property is known as \emph{robustness} of the Bayes procedure.

This paper summarizes recent results \cite{OwhadiScovel:2013, OSS:2013} that give conditions under which Bayesian inference appears to be non-robust in the most extreme fashion, in the sense that arbitrarily small changes of the prior and model class lead to arbitrarily large changes of the posterior value of a quantity of interest.
We call this extreme non-robustness ``brittleness'', and it can be visualized as the smooth dependence of the value of the quantity of interest on the prior breaking into a fine patchwork, in which nearby priors are associated to diametrically opposed \emph{posterior} values.
Naturally, the notion of ``nearby'' plays an important role, and this point will be revisited later.

Much as classical numerical analysis shows that there are `stable' and `unstable' ways to discretize a partial differential equation, these results and the wider literature of positive \cite{Bernstein:1964, CastilloNickl:2013, Doob:1949, Kleijnvanderv:2012, LeCam:1953, Stuart:2010, vonMises:1964} and negative results \cite{Belot:2013, DiaconisFreedman:1986, Freedman:1963, Freedman:1999, Johnstone:2010, Leahu:2011} on Bayesian inference contribute to an emerging understanding of `stable' and `unstable' ways to apply Bayes' rule in practice.

The results reported in this paper show that the process of Bayesian conditioning on data at finite enough resolution is unstable (or ``sensitive'' as defined in \cite{TibshiraniWasserman:1988}) with respect to the underlying distributions (under the total variation
and Prokhorov metrics) and the source of negative results similar to those caused by tail properties in statistics
\cite{BahadurSavage:1956, Donoho:1988}. The mechanisms causing the stability/instability of posterior predictions suggest that learning and robustness are conflicting requirements and raise the question of a missing stability condition for using Bayesian inference for continuous systems with finite information  (akin to the CFL stability condition for using discrete schemes to approximate continuous PDEs).

\subsection*{Bayes' Theorem and Robustness}

To begin, consider a simple example of Bayesian reasoning in action:
\begin{pb}
	\label{pb:1}
	There is a bag containing 102 coins, one of which always lands on heads, while the other 101 are perfectly fair.  One coin is picked uniformly at random from the bag, flipped 10 times, and 10 heads are obtained.  What is the probability that this coin is the unfair coin?
\end{pb}

The correct probability is given by applying Bayes' theorem:
\begin{equation}
	\label{eq:Bayes}
	\P[A|B]=\P[B|A]\frac{\P[A]}{\P[B]}=\frac{1}{1+101 \times 2^{-10}} \approx 0.91,
\end{equation}
where $A$ is the event `the coin is the unfair coin' and $B$ is the event `10 heads are observed'.
If the number of coins is not known exactly and the supposedly fair coins are not exactly fair, then Bayes' theorem   produces a robust inference in the following sense:  if the fair coins are slightly unbalanced and the probability of getting a tail is  $0.51$, and an estimate of 100 coins is used and an estimate $\frac{1}{2}$ of the fairness of the fair coins is used, then the resulting estimate $\frac{1}{1+99 \times 2^{-10}}$ is still a good approximation to the correct answer.
Observe also that if the prior estimate of the number of coins in the bag is grossly wrong (e.g.~$10^6$) then the posterior would still be accurate in the limit of infinitely many coin flips:  in this case, the Bayesian estimator is said to be \emph{consistent}.

Do these conclusions remain true when the underlying probability space is continuous or an approximation thereof?
For example, what if the random outcomes are decimal numbers --- perhaps given to finite precision --- rather than heads or tails?

\subsection*{The General Problem and its Bayesian Answer}
\label{sec_general}

\begin{pb}
	\label{pb:2}
	Let $\mathcal{X}$ denote the space in which observations/samples take their values, and let $\mathcal{M}(\mathcal{X})$ denote the set of probability measures on $\mathcal{X}$.  Let $\Phi \colon \mathcal{M}(\mathcal{X})\to \R$ be a function\footnote{All spaces will be topological spaces, the term `function' will mean Borel measurable function and `measure' will mean Borel measure.} defining a \emph{quantity of interest}.  Let the \emph{data-generating distribution}  $\mu^\dagger \in \mathcal{M}(\mathcal{X})$ be an unknown or partially known probability measure on $\mathcal{X}$.  The objective is to estimate $\Phi(\mu^\dagger)$ from the observation of $n$ i.i.d.~samples from $\mu^\dagger$, which we denote by $d=(d_1,\ldots,d_n)\in \mathcal{X}^n$.
\end{pb}

\begin{eg}
	When $\mathcal{X}$ is the real line $\R$, a prototypical example of a quantity of interest is $\Phi(\mu):=\mu[X\geq a]$, the probability that the random variable $X$ distributed according to $\mu$ exceeds the threshold value $a$.
	However, the results that we report below apply to \emph{any} pre-specified quantity of interest $\Phi$.
\end{eg}

The Bayesian answer to this problem is to model $\mu^\dagger$'s generation of sample data as coming from a random measure on $\mathcal{X}$ and to condition $\Phi$ with respect to the observation of the $n$ i.i.d.\ samples.
This is done by choosing a \emph{model class} $\mathcal{A}\subseteq\mathcal{M}(\mathcal{X})$  and a probability measure $\pi \in \mathcal{M}(\mathcal{A})$ which we call \emph{the prior}.
This prior determines the randomness with which a representative $\mu \in \mathcal{A}$ is selected, and, for each such $\mu \in \mathcal{A}$, the generation of $n$ i.i.d.\ samples $d \in \mathcal{X}^n$ by randomly sampling from $\mu^{n}$ naturally determines a product measure on $\mathcal{A} \times \mathcal{X}^n$.
The prior estimate of the quantity of interest is $\E_{\mu\sim \pi} [\Phi(\mu)]$ and, for an open\footnote{We assume $B$ to be open and of strictly positive measure to avoid problems associated with conditioning with respect to events of measure zero.} $B \subseteq \mathcal{X}^n$, the posterior estimate is defined as the conditional expectation $\E_{\mu \sim \pi, d\sim \mu^{n}} [\Phi(\mu)| d\in B]$ with respect to this product measure.

The connection to the standard presentation of Bayesian inference in terms of a prior on a parameter space is as follows:
to construct a model class $\mathcal{A} \subseteq \mathcal{M}(\mathcal{X})$ and a prior $\pi_0 \in \mathcal{M}(\mathcal{A})$ from a Bayesian parametric model $\mathcal{P} \colon \Theta\to \mathcal{M}(\mathcal{X})$ defined on a \emph{parameter space} $\Theta$ equipped with a prior $p_{0}\in \mathcal{M}(\Theta)$, one simply pushes forward under the map $\mathcal{P}$.
That is, the model class $\mathcal{A} \subseteq \mathcal{M}(\mathcal{X})$ is defined by $\mathcal{A}:=\mathcal{P}(\Theta)$ and the prior $\pi_0 \in\mathcal{M}(\mathcal{A})$ is defined as the push-forward $\pi_0:=\mathcal{P}p_{0}$ of $p_0$ by the model $\mathcal{P}$, i.e.\ $\pi_{0}(E):= p_{0}(\mathcal{P}^{-1}(E))$ for measurable $E \subseteq \mathcal{A}$.

\subsection*{Inconsistency under Mis-specification}

We now discuss the effects of mis-specification for a Bayesian parametric model $\mathcal{P} \colon \Theta \to \mathcal{M}(\mathcal{X})$.
It is convenient to denote such a model by $\mathcal{P} \colon \theta \mapsto \mu(\theta)$, so that the model class is $\mathcal{A} := \mathcal{P}(\Theta) = \{\mu(\theta) \mid \theta \in\Theta\}$.
If the model class $\mathcal{P}(\Theta)$ contains the data-generating distribution $\mu^\dagger$, i.e.~if there is some parameter value $\theta \in \Theta$ such that $\mu^\dagger = \mu(\theta)$, then the model is said to be \emph{well-specified};  otherwise, it is said to be \emph{mis-specified}.

For simplicity, consider the classical case where, for each $\theta \in\Theta$, $\mu(\theta)$ has a probability density function with respect to some common reference measure on $\mathcal{X}$, that is $\mu(\theta)=p(\quark,\theta) \, \mathrm{d}x$ for some measure $\mathrm{d}x$.
Then, for a prior $p_{0}\in \mathcal{M}(\Theta)$, let $p_n \in \mathcal{M}(\Theta)$ denote the posterior distribution on $\Theta$ after observing the data $d$ (see e.g.~\cite[p.~126]{Berger:1985}) and push forward both the prior and posterior to their corresponding measures, $\pi_0:=\mathcal{P}p_0$ and $\pi_n:=\mathcal{P}p_n$, on $\mathcal{M}(\mathcal{A})$.

Now suppose that the model is well-specified  and that $p_0$ gives strictly positive mass to every neighborhood of every point $\theta \in \Theta$ --- this assumption of `maximal open-mindedness' is commonly referred to as \emph{Cromwell's rule} \cite{Lindley:1985}.
Then, when $\Theta$ is finite dimensional, under suitable regularity conditions, the posterior value of the quantity of interest $\E_{\mu \sim \pi_n}\bigl[\Phi(\mu)\bigr]$ converges to $\Phi(\mu^\dagger)$  as $n\to \infty$.
This convergence, which can be shown to be asymptotically normal, is commonly referred to as the \emph{Bernstein--von Mises theorem} or \emph{Bayesian central limit theorem} \cite{Bernstein:1964, Doob:1949, LeCam:1953, vonMises:1964}.
However, for infinite-dimensional $\Theta$ and with similar regularity and strict positivity assumptions, there is a wealth of positive \cite{CastilloNickl:2013, Kleijnvanderv:2012, Stuart:2010} and negative results \cite{Belot:2013, DiaconisFreedman:1986, Freedman:1963, Freedman:1999, Johnstone:2010, Leahu:2011} showing that the truth/falsity of the Bernstein--von Mises property depends sensitively on subtle topological and geometrical details.

Conversely, if the model is mis-specified, then, under regularity conditions \cite{Berk:1966, KleijnVaart:2006, Kleijnvanderv:2012, Shalizi:2009}, the posterior value $\E_{\mu \sim \pi_n}\bigl[\Phi(\mu)\bigr]$ converges as $n \to \infty$ to $\Phi\big(\mu(\theta^{\ast})\big)$, where $\theta^{\ast}$ maximizes the expected log-likelihood function $\theta \mapsto \E_{\mu^{\dagger}} \bigl[ \log p(\quark,\theta) \bigr]$.
If, in addition,  $\mu^\dagger$ is absolutely continuous with respect to each $\mu(\theta)$ for $\theta \in \Theta$, then $\theta^{\ast}$ can also be shown to minimize the \emph{Kullback--Leibler divergence} or \emph{relative entropy distance} $\theta \mapsto D_{\mathrm{KL}}\bigl( \mu^{\dagger} \big\| \mu(\theta) \bigr)$ from $\mu^\dagger$ to $\mu(\theta)$.

\begin{eg}
	\label{eg:Gaussian_tail}
	To illustrate this, let $\mathcal{X}=\R$ and consider the Gaussian model where $\mu(c,\sigma)$ is a Gaussian with mean $c$ and standard deviation $\sigma$, that is, it has the probability density
	\[
		p(x,c,\sigma):=\frac{1}{\sigma \sqrt{2\pi}}\exp \left( -\frac{(x-c)^2}{2 \sigma^2}\right),
	\]
	and therefore the expected log-likelihood is
	\[
		\E_{\mu^{\dagger}} \bigl[ \log p(\cdot,c,\sigma) \bigr] = - \int_{\R}  \frac{(x - c)^{2}}{2 \sigma^{2}} \, \mathrm{d} \mu^{\dagger}(x) - \log \sigma  - \log \sqrt{2 \pi}.
	\]
	If, for a data-generating distribution $\mu^{\dagger}$ with finite second moments, we let $c^{\dagger}$ denote its mean and $\sigma^{\dagger}$ its standard deviation, then a quick calculation shows that $\theta^{\ast} = (c^{\ast}, \sigma^{\ast})$ maximizes the expected log-likelihood if and only if $c^{\ast} = c^{\dagger}$ and $\sigma^{\ast} = \sigma^{\dagger}$.
	Hence, the asymptotic Bayesian posterior estimate of $\Phi(\mu^\dagger)$ is $\Phi\big(\mu(c^\dagger,\sigma^\dagger)\big)$, irrespective of what the quantity of interest $\Phi$ might be.
	However, there are \emph{many} different probability distributions $\mu$ on $\R$ that have the same first and second moments as $\mu^{\dagger}$ but have different higher-order moments, or different quantiles.
	Predictions of those other moments or quantiles using the Gaussian distribution $\mu (c^{\dagger}, \sigma^{\dagger})$ can be inaccurate by orders of magnitude.
	A simple example is provided by the tail probability $\Phi(\mu) := \P_\mu \bigl[ |X-c_\mu| \geq t \sigma_\mu \bigr]$, where $c_\mu$ and $\sigma_\mu$ denote the mean and standard deviation of $\mu$, and $t > 0$.
	Under the Gaussian model,
	\[
		\P_\mu \bigl[ | X - c_\mu | \geq t \sigma_\mu \bigr] = 1 + \mathrm{erf} \left( - \frac{t}{\sqrt{2}} \right),
	\]
	whereas the extreme cases that prove the sharpness of Chebyshev's inequality --- in which the probability measure is a discrete measure with support on at most three points in $\R$ --- have
	\[
		\P_\mu \bigl[ | X - c_\mu | \geq t \sigma_\mu \bigr] = \min \left\{ 1, \frac{1}{t^{2}} \right\}.
	\]
	In the case of the archetypically rare `$6 \sigma$ event', i.e.~$t = 6$, the ratio between the two is approximately $1.4 \times 10^{7}$.
	This comparison is, of course, almost perversely extreme:  it would be obvious to any observer with only moderate amounts of `Chebyshev-type' sample data that the data were being drawn from a highly non-Gaussian distribution.
	However, it is not inconceivable that the true distribution $\mu^{\dagger}$ has a Gaussian-looking bulk but also has tails that are significantly fatter than those of a Gaussian, and the difference may be difficult to establish using reasonable amounts of sample data;  yet, it is those tails that drive the occurrence of `Black Swans', catastrophically high-impact but low-probability outcomes.
\end{eg}

Although it is understood that Bayesian estimators can be inconsistent if the model is grossly mis-specified, a pressing question is whether they have good convergence properties when the model class  $\{ \mu(\theta) \mid \theta\in \Theta\}$ is `close enough' to the truth  $\mu^{\dagger}$ in an appropriate sense.

Such concerns can be traced back to Box's dictum that ``essentially, all models are wrong, but some are useful'' \cite[p.~424]{Box:1987} and question ``how wrong do they have to be to not be useful?'' \cite[p.~74]{Box:1987}.
These queries are also critical because, although gross mis-specification of the model can be detected before engaging in a complete Bayesian analysis \cite{HausmanTaylor:1981, White:1982}, usually one \emph{cannot be sure} that the model is well-specified.

To answer these questions we will examine the robustness of Bayesian inference by computing optimal bounds on prior and posterior values in terms of given sets of priors.
Indeed, the exploration of classes of Bayesian models is one response to the concern that the choice of prior-likelihood combination could, to some degree, be arbitrary, and forms the basis of the approach known as \emph{robust Bayesian inference} \cite{Berger:1984, Berger:1994, Box:1953, WassermanEtAl:1993, Wasserman:1990}.
To do so, we need some definitions.

\begin{defn}
	\label{defn:robustness}
	For a model class $\mathcal{A}\subseteq \mathcal{M}(\mathcal{X})$, a quantity of interest $\Phi \colon \mathcal{A} \to \R$, and a set of priors $\Pi \subseteq \mathcal{M}(\mathcal{A})$, let
	\begin{align*}
		\mathcal{L}(\Pi)&:=\inf_{\pi \in \Pi}\E_{\mu\sim \pi}\bigl[\Phi(\mu)\bigr]  \\
  		\mathcal{U}(\Pi)&:=\sup_{\pi \in \Pi}\E_{\mu\sim \pi}\bigl[\Phi(\mu)\bigr]
	\end{align*}
	denote the \emph{optimal lower and upper bounds on the prior values} of $\Phi$.  For $B$ a non-empty open subset of the data space $\mathcal{X}^n$, let $\Pi_B\subseteq \Pi$ be the subset of priors $\pi$ such that the probability that $d\in B$ is nonzero, i.e.~$\mathbb{P}_{\mu \sim \pi, d\sim \mu^{n}}[d \in B]>0$, and let
	\begin{align*}
		\mathcal{L}(\Pi|B)&:=\inf_{\pi\in \Pi_B }\E_{\mu \sim \pi, d\sim \mu^{n}} \bigl[\Phi(\mu)\big| d\in B\bigr] \\
		\mathcal{U}(\Pi|B)&:=\sup_{\pi\in \Pi_B }\E_{\mu \sim \pi, d\sim \mu^{n}} \bigl[\Phi(\mu)\big| d\in B\bigr]
	\end{align*}
	denote the \emph{optimal lower and upper bounds on the posterior values} of $\Phi$ given that $d\in B$.
\end{defn}

\subsection*{Brittleness under Infinitesimal Perturbations}

\begin{figure}[t!]
	\begin{center}
		\includegraphics[width=0.6\textwidth]{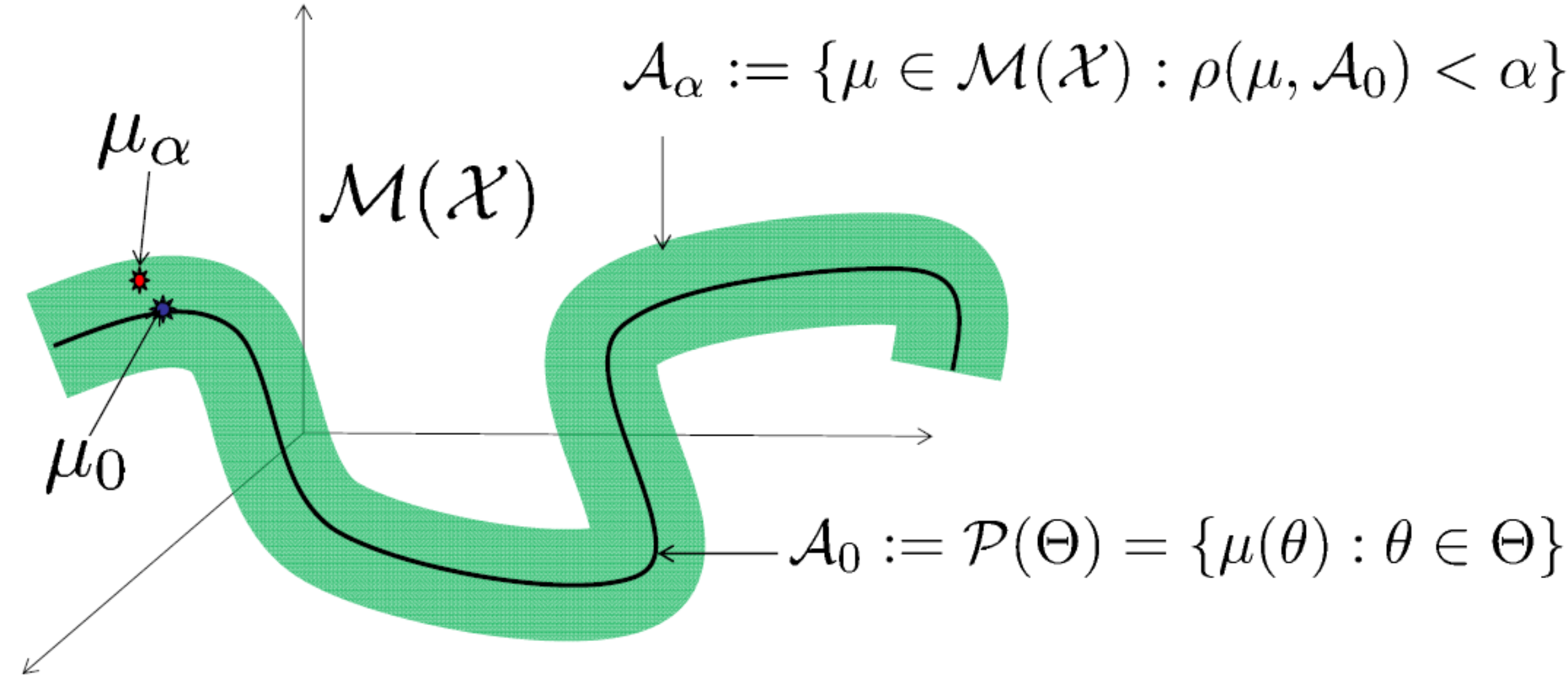}
	\end{center}
	\caption{The original model class $\mathcal{A}_{0}$ (black curve) is enlarged to its metric neighbourhood $\mathcal{A}_{\alpha}$ (shaded).  This procedure determines perturbations $\mu_{\alpha} \in \mathcal{A}_{\alpha}$ of the original random measure $\mu_{0} \in \mathcal{A}_{0}$.}
	\label{fig:perturbation}
\end{figure}

Consider again the  model $\mathcal{P} \colon \Theta\to \mathcal{M}(\mathcal{X})$ but now denote the model class by $\mathcal{A}_{0} := \mathcal{P}(\Theta)=\{\mu(\theta) \mid \theta \in \Theta\}$ and the prior by
$\pi_0 \in  \mathcal{M}\bigl(\mathcal{A}_{0}\bigr)$.
To quantify perturbations in the model, and define what it means for two distributions to be close to one another, we select a metric $\rho$ on $\mathcal{M}(\mathcal{X})$.
As illustrated in Figure \ref{fig:perturbation}, for $\alpha > 0$, we enlarge the set $\mathcal{A}_{0}$ to its metric neighborhood $\mathcal{A}_{\alpha}$ and thereby naturally determine a set of priors $\Pi_{\alpha} \subseteq \mathcal{M}(\mathcal{A}_{\alpha})$ such that the random measure $\mu_{\alpha}$ associated with every $\pi_{\alpha} \in \Pi_{\alpha}$ lies within distance $\alpha$ from the random measure $\mu_{0}$ associated with the prior $\mu_{0}$ and the Bayesian model $\mathcal{P}$.
Then we analyze the robustness of its posteriors, as in Definition \ref{defn:robustness},  with respect to these size-$\alpha$ perturbations.

To that end, suppose that $\mathcal{X}$ is metrizable and select a consistent metric $\mathrm{d}$ for $\mathcal{X}$.  Let $\mathcal{B}(\mathcal{X})$ denote the Borel subsets of $\mathcal{X}$.  We will consider two metric distances $\rho(\mu, \nu)$ between $\mu, \nu \in \mathcal{M}(\mathcal{X})$:  $\rho$ will be either the \emph{total variation (TV) metric}
\[
	\rho_{\textup{TV}}(\mu,\nu):=\sup \bigl\{ | \mu(A)-\nu(A) | \big| A \in \mathcal{B}(\mathcal{X})   \bigr\},
\]
or the \emph{Prokhorov metric}\footnote{The TV metric is generally considered to generate too strong a topology on the space $\mathcal{M}(\mathcal{X})$ of probability measures, and the weak topology is generally considered more appropriate, see e.g.~\cite{Billingsley1}.
Fortunately, when $\mathcal{X}$ is separable, this topology is metrized by the Prokhorov metric. For a thorough discussion regarding metrics on spaces of measures see e.g.~\cite{Rachev}.}
\[
	\rho_{\textup{P}}(\mu,\nu):=\inf \left\{ \varepsilon > 0 \,\middle|\,
		\mu(A) \leq \nu(A^\varepsilon) + \varepsilon,\, A\in \mathcal{B}(\mathcal{X})  \right\},
\]
where $A^\varepsilon := \{x \in \mathcal{X}| \mathrm{d}(x,A) < \varepsilon \}$.
For $\alpha>0$, the neighborhood $\mathcal{A}_{\alpha}$ of $\mathcal{A}_{0}$ emerges naturally from the ball fibration
\[
	\mathcal{A}^{\ast}:= \bigl\{ (\mu_{1},\mu_{2}) \in \mathcal{M}(\mathcal{X}) \times \mathcal{M}(\mathcal{X}) \,\big|\, \mu_{1}\in \mathcal{A}_{0},  \rho(\mu_{2},\mu_{1})<\alpha \bigr\},
\]
in the sense that if $P_0$ and $P_\alpha$ denote the projections onto the first and second components of $\mathcal{M}(\mathcal{X}) \times \mathcal{M}(\mathcal{X})$, then $P_{0}\mathcal{A}^{\ast}=\mathcal{A}_{0}$ and  $P_{\alpha}\mathcal{A}^{\ast}=\mathcal{A}_{\alpha}$.
Consequently, a natural set of priors $\Pi_{\alpha} \subseteq \mathcal{M}(\mathcal{A}_{\alpha})$ corresponding to $\pi_{0} \in \mathcal{M}(\mathcal{A}_{0})$ is defined by
\[
	\Pi_{\alpha}:= \bigl\{ \pi_\alpha \in \mathcal{M}(\mathcal{A}_\alpha) \,\big|\, \text{for some $\pi \in \mathcal{M}(\mathcal{A}^{\ast})$, $P_0\pi= \pi_{0}$ and $P_\alpha\pi=\pi_\alpha$} \bigr\}.
\]

To state our result, consider again Problem \ref{pb:2}, and let some $x^n:=(x_{1},\dots,x_{n}) \in \mathcal{X}^{n}$ be a point such that we observe that $d\in B^{n}_{\delta}:=\prod_{i=1}^{n}{B_{\delta}(x_{i})}$, where $B_{\delta}(x)\subseteq \mathcal{X}$ is the open ball of radius $\delta$ centered on $x \in \mathcal{X}$.
Using the notation of Definition \ref{defn:robustness}, and $\Pi_{\alpha}$ defined above in terms of the  TV  or Prokhorov metric, the Brittleness Theorem 6.4 of \cite{OSS:2013} then reads as follows\footnote{All results of this paper and those in  \cite{OwhadiScovel:2013, OSS:2013, OSSMO:2011} require some mild technical measure-theoretic and topological assumptions.
For example, here it is sufficient if $\mathcal{P}(\Theta)$ is a Borel subset of a  Polish space (a separable completely metrizable space).
Unfortunately, $\mathcal{M}(\mathcal{X})$ is not generally separable with respect to the TV metric, and hence is not Polish.
However, if $\mathcal{X}$ is Polish, then $\mathcal{M}(\mathcal{X})$ topologized by weak convergence is Polish and the Prokhorov metric provides a complete metrization of it.
Consequently, when $\Theta$ is Polish, $\mathcal{X}$ is Polish, and $\mathcal{P}$ is injective and measurable with respect to the weak topology, it then follows from Suslin's Theorem that $\mathcal{P}(\Theta)$ is a Borel subset of the Polish space $\mathcal{M}(\mathcal{X})$.
For a thorough investigation of such matters, illustrating the benefits of Polish spaces as the foundation for the framework, see \cite{OSS:2013}.}:

\begin{thm}\label{thm:1}
	If
	\begin{equation}
		\label{eq:zeromass}
		\lim_{\delta \downarrow 0} \sup_{x \in \mathcal{X}} \sup_{\theta\in \Theta} \mu(\theta)[B_{\delta}(x)]=0,
	\end{equation}
	then, for all $\alpha>0$, there exists $\delta(\alpha)>0$ such that for all $0<\delta<\delta(\alpha)$, all $n \in \mathbb{N}$ and all $x^n \in \mathcal{X}^n$,
	\[
		\mathcal{L}(\Pi_{\alpha}|B^{n}_{\delta})\leq \operatorname{ess\,inf}_{\pi_0} (\Phi)\quad\text{ and } \quad\operatorname{ess\,sup}_{\pi_0} (\Phi) \leq  \mathcal{U}(\Pi_{\alpha}|B^{n}_{\delta})\, ,
	\]
	where  $\operatorname{ess\,inf}_{\pi_0} (\Phi):=\sup\{r\mid \pi_0 [\Phi<r]=0\}$ and $\operatorname{ess\,sup}_{\pi_0} (\Phi):=\inf\{r\mid \pi_0 [\Phi>r]=0\}$.
\end{thm}

\smallskip

Note that  condition \eqref{eq:zeromass} is extremely weak and satisfied for most parametric Bayesian models.
Furthermore, suppose that Cromwell's rule is applied.
Then, although it implies consistency if the model is well-specified, here it leads to maximal brittleness under local mis-specification.
More precisely, under Cromwell's rule, $\operatorname{ess\,inf}_{\pi_0} (\Phi)=\inf_{\mu \in \mathcal{A}_0} \Phi(\mu)$ and $\operatorname{ess\,sup}_{\pi_0} (\Phi)=\sup_{\mu \in \mathcal{A}_0} \Phi(\mu)$, so the conclusion of Theorem \ref{thm:1} becomes
\[
	\mathcal{L}(\Pi_{\alpha}|B^{n}_{\delta})\leq \inf_{\mu \in \mathcal{A}_0} \Phi(\mu)
\quad \text{ and } \quad  \sup_{\mu \in \mathcal{A}_0} \Phi(\mu) \leq \mathcal{U}(\Pi_{\alpha}|B^{n}_{\delta})\, .
\]
In other words, the range of posterior predictions among all admissible priors is as wide as the deterministic range of the quantity of interest $\Phi$.

Note that since $\Phi$ is arbitrary, the brittleness described in Theorem \ref{thm:1} is not limited to a quantile or moment of $\mu$ but concerns its whole posterior distribution.

\subsection*{Brittleness under Finite Information}
\label{con:psi}
One response to the concern  that the choice of prior and model are somewhat arbitrary
 \cite{WassermanEtAl:1993} is  to perform a sensitivity analysis over classes of priors and models.
One way to specify a class $\Pi$ of admissible priors $\pi$ is to select some `features' (such as the polynomial moments, or other functionals) and specify some values, ranges, or distributions for those features.
It is interesting to understand the impact of the features left unspecified, i.e.~the \emph{codimension} and not just the \emph{dimension} of $\Pi$;
while \emph{robust Bayesian inference} \cite{Berger:1984, Berger:1994, Box:1953, Wasserman:1990} has shown that posterior conclusions remain stable when $\Pi$ is finite-dimensional, our results can be interpreted as saying that brittleness ensues whenever $\Pi$ has finite codimension, regardless of how large its codimension is. It is important to note that this is in some sense the generic situation:  when $\mathcal{A}$ is an infinite set, one would have to specify infinitely many features of priors $\pi \in \Pi$ to achieve a finite-dimensional $\Pi$;
from a computational and epistemic standpoint, the specification of infinitely many features in finite time appears to be somewhat problematic.

To study this problem, we introduce a representation space $\mathcal{Q}$ (e.g., prototypically, $\R^{k}$) and a mapping $\Psi \colon \mathcal{A}\to \mathcal{Q}$ from the subset $\mathcal{A} \subseteq \mathcal{M}(\mathcal{X})$ into $\mathcal{Q}$, which can be thought of as a map to `generalized moments'.
Let $\mathfrak{Q} \subseteq \mathcal{M}(\mathcal{Q})$ be a subset of the set of probability distributions on $\mathcal{Q}$ such that each distribution $\mathbb{Q}\in \mathfrak{Q}$ has its support contained in $\Psi(\mathcal{A})$.
If the set $\mathfrak{Q}$ represents priors for the distribution of $\Psi(\mu), \mu \in \mathcal{A}$, then a naturally induced set of priors $\Pi$ on $\mathcal{A}$ is the pull-back $\Pi:=\Psi^{-1}(\mathfrak{Q}) \subseteq \mathcal{M}(\mathcal{A})$, defined by $\Psi^{-1}(\mathfrak{Q}) :=\{\pi \in \mathcal{M}(\mathcal{A}) \mid \Psi \pi \in \mathfrak{Q}\}$.

\begin{eg}
	\label{eg:1}
	Consider the case $\mathcal{X}=[0,1]$, $\mathcal{A}:=\mathcal{M}([0,1])$, and $\Phi(\mu)=\E_{\mu}[X]$.
	Thus, the aim is to estimate the mean $\Phi(\mu^{\dagger})=\E_{\mu^{\dagger}}[X]$ of the random variable $X$ corresponding to some unknown measure $\mu^\dagger\in \mathcal{A}$ and we observe $d=(d_1,\ldots,d_n)$, $n$ i.i.d.\ samples from $X$.
	Let $k$ be fixed and let $\Psi(\mu)=(\E_{\mu}[X],\ldots,\E_{\mu}[X^{k}])$ be the map to the first $k$ polynomial moments.
	If we write a point $q\in \R^{k}$ in terms of its coordinates $q:=(q_1,\ldots,q_{k})$, then $\Psi^{-1}(q)$ is exactly the set of measures $\mu \in \mathcal{M}([0,1])$ such that $\E_{\mu}[X^i]=q_i$ for $1 \leq i \leq k$.
	Now define a measure $\mathbb{Q}$ on the truncated moment space $\Psi(\mathcal{M}([0,1]) \subseteq \R^{k}$ as follows.
	Since the first moment $\E_{\mu}[X], \mu \in \mathcal{M}([0,1])$, ranges over the unit interval, consider the uniform measure on the unit interval in the first coordinate.
	Next define the conditional measure when the first coordinate is $q_1 \in [0,1]$ to be uniform on the range of the second moment $\left[ \inf_{\mu:\, \E_{\mu}[X]=q_1} \E_\mu[X^2], \sup_{\mu:\, \E_{\mu}[X]=q_1} \E_\mu[X^2] \right]$.
	Repeat this conditioning process on the higher coordinates iteratively in the same manner.
	Then, the induced set of priors $\Pi:=\Psi^{-1}\mathbb{Q}$ on $\mathcal{M}([0,1])$ is the set of measures $\pi$ such that, when $\mu \sim \pi$, the distribution of $(\E_{\mu}[X],\ldots,\E_{\mu}[X^{k}])$ is $\mathbb{Q}$.
\end{eg}

We now state the Brittleness Theorem 4.13 in \cite{OSS:2013} for the general case of Problem \ref{pb:2},
and apply it to Example \ref{eg:1}.
To that end, let the model class $\mathcal{A}\subseteq \mathcal{M}(\mathcal{X})$ be chosen along with a generalized moment map $\Psi \colon \mathcal{A}\to \mathcal{Q}$ to a representation space $\mathcal{Q}$.
Let $\mathfrak{Q} \subseteq \mathcal{M}(\mathcal{Q})$ be a specified set of priors on $\mathcal{Q}$ and from them determine  $\Pi:=\Psi^{-1}(\mathfrak{Q})\subseteq \mathcal{M}(\mathcal{A})$ as the induced set of priors. For fixed $(x_{1},\dots,x_{n}) \in \mathcal{X}^{n}$, let $B^{n}_{\delta}:=\prod_{i=1}^{n}{B_{\delta}(x_{i})}$ where $B_{\delta}(x)$ is the open ball of radius $\delta$ centered on $x \in \mathcal{X}$.
The following theorem gives optimal bounds on posterior values for the class of priors $\Pi$ defined above, given that the observation $d\in B^{n}_{\delta}$.

\begin{thm}
	\label{thm:2}
	Suppose that, for all $\gamma >0$, there exists some $\mathbb{Q}\in \mathfrak{Q}$ such that
	\begin{equation}
		\label{eq:cond1}
		\E_{q\sim \mathbb{Q}} \left[ \inf_{\mu\in  \Psi^{-1}(q), i=1,\dots,n} \mu[B_{\delta}(x_{i})] \right]=0
	\end{equation}
	and
	\begin{equation}
		\label{eq:cond2}
		\mathbb{P}_{q\sim \mathbb{Q}} \left[ \sup_{\mu \in \Psi^{-1}(q):\, \mu[B_{\delta}(x_{i})]>0, i=1,\dots,n}\Phi(\mu) > \sup_{\mu\in \mathcal{A}}\Phi(\mu) - \gamma \right]>0\, .
	\end{equation}
	Then
	\[
		\mathcal{U}\big(\Pi\big|B^{n}_{\delta}\big) =\sup_{\mu \in \mathcal{A}}\Phi(\mu),
	\]
	with similar expressions for the lower bounds $\mathcal{L}$.
\end{thm}

In other words, if there is a measure $\mathbb{Q} \in \mathfrak{Q}$ such that for  $\mathbb{Q}$-almost all $q \in \mathcal{Q}$ there is a $\mu \in \Psi^{-1}(q)$ which achieves an arbitrarily small mass on one of $B_{\delta}(x_{i}), i=1,\dots,n$, and  with  non-zero $\mathbb{Q}$ probability there is $\mu \in \Psi^{-1}(q)$ which almost extremizes $\Phi$ while putting positive mass on all $B_{\delta}(x_{i}), i=1,\dots,n$, then the range $\bigl[ \mathcal{L}\big(\Pi\big|B^{n}_\delta\big) , \mathcal{U}\big(\Pi\big|B^{n}_\delta\big) \bigr]$ of posterior values for $\Phi$ is exactly the `deterministic' range of $\Phi$, i.e.~$\bigl[ \inf_{\mu\in \mathcal{A}}\Phi(\mu) , \sup_{\mu\in \mathcal{A}}\Phi(\mu) \bigr]$.

Conditions \eqref{eq:cond1} and \eqref{eq:cond2} are very weak and simple dimensionality arguments suggest that they are typically satisfied if $\mathcal{Q}$ is finite-dimensional.
Hence, although Bayesian inference is robust in situations where  the distributions of \emph{all but} finitely many generalized moments of the data-generating distribution $\mu^{\dagger}$ are known, Theorem \ref{thm:2} suggests that it is  brittle when the distributions of \emph{only} finitely many generalized moments of $\mu^{\dagger}$ are known, while infinitely many remain unknown.
As an example, it is instructive to observe how Theorem \ref{thm:2}, applied to Example \ref{eg:1} in \cite[Ex.~4.16]{OSS:2013}, shows that if the data-generating measure has some non-atomic component, then when the number of  samples $n$ is large enough and $\delta$ small enough, the optimal bounds on posterior values of $\Phi(\mu)=\E_{\mu}[X]$, given the distribution $\mathbb{Q}$ defined on its moments, are $0$ and $1$.

To quantify ``large enough'' and ``small enough'' and to remove the ``non-atomic'' requirement above, Theorem~3.1 of \cite{OwhadiScovel:2013} provides a quantitative version of Theorem~\ref{thm:2} in which the conditions of the theorem are only required to hold approximately.
When applied to Example \ref{eg:1} with the set $\Pi:=\Psi^{-1}\mathbb{Q}$ of priors generated instead by the uniform prior $\mathbb{Q}$ restricted to the truncated moment space, Theorem 3.3 of \cite{OwhadiScovel:2013}  establishes that, although the prior value satisfies $\mathcal{U}(\Pi)=\frac{1}{2}$, the posterior value satisfies
\begin{equation}\label{eqkjhkdhjdh}
	1 - 4 e \Bigl( \frac{2k \delta}{e} \Bigr)^\frac{1}{2k+1} \leq \mathcal{U} \big(\Pi|B_{\delta}\big) \leq 1\, .
\end{equation}
Consequently, regardless of the number of moment constraints $k$, and the location of a single data point, for $\delta$ smaller than an elementary known function of $k$, we have brittleness. This result also holds for arbitrary multiple samples.
Remark 4.18 of \cite{OSS:2013} also suggests that brittleness would persist if the hard bound $\delta$ to specify measurement uncertainty is replaced by a level of noise with variance decreasing with $\delta$.

\subsection*{Mechanism Causing Brittleness}

We will now illustrate one mechanism causing brittleness with a simple example derived from the proof of
Theorem \ref{thm:1}.
In this example we are interested in estimating $\Phi(\mu^\dagger)=\E_{\mu^\dagger}[X]$ where $\mu^\dagger$ is an unknown distribution on the unit interval ($\mathcal{X}=[0,1]$)  based on the observation of a single data point $d_1=\frac{1}{2}$ up to resolution $\delta$ (i.e.\ we observe $d_1\in B_\delta(x_1)$ with $x_1=\frac{1}{2}$).

\begin{figure}[tp]
	\begin{center}
		\subfigure[$f^a(x,\theta)$]{
			\includegraphics[height=2.5cm]{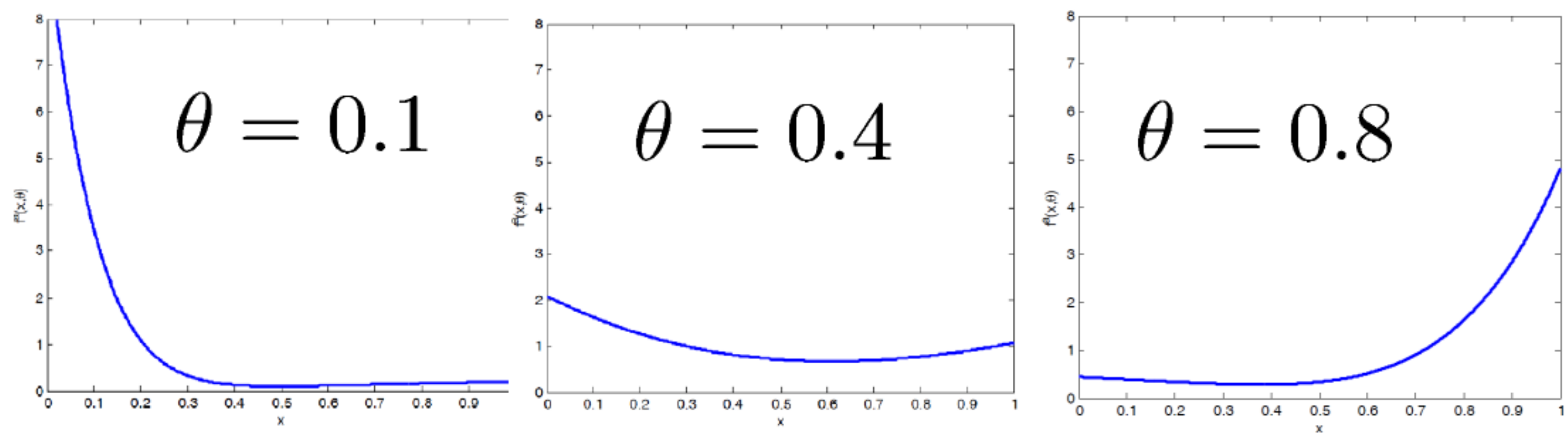}
		}
		\subfigure[$f^b(x,\theta)$]{
			\includegraphics[height=2.4cm]{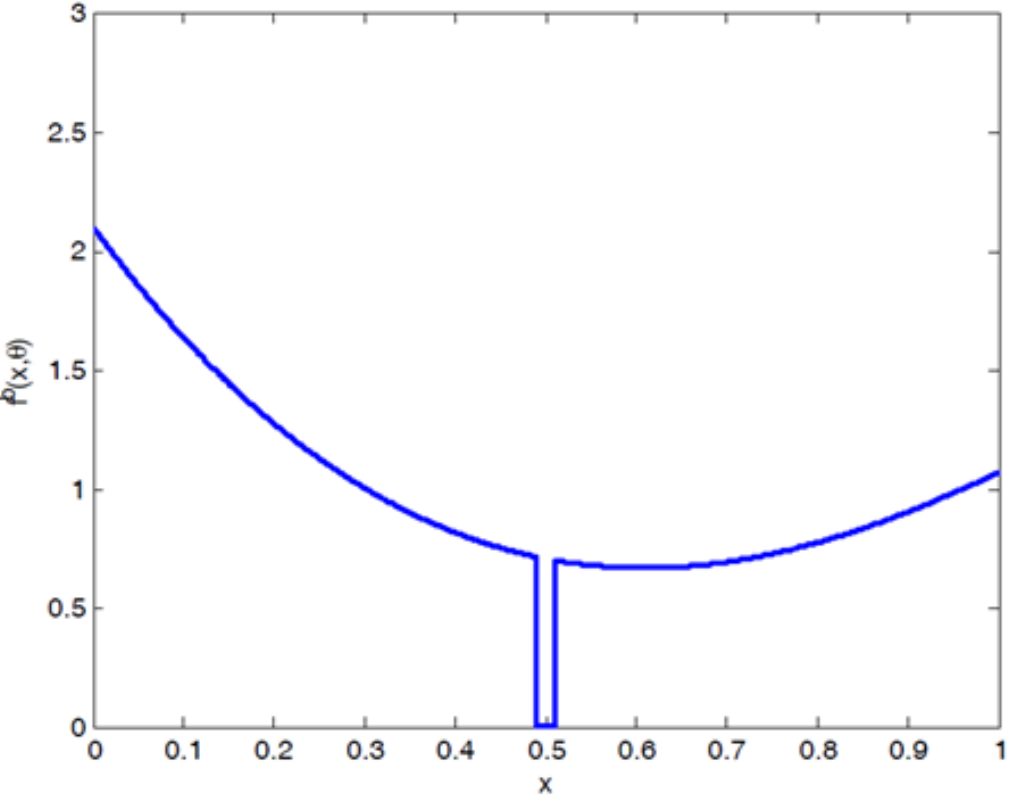}
		}
	\end{center}
	\caption{Illustration of the density $f^a(x,\theta)$ of model $a$ and $f^b(x,\theta)$ of model $b$.}
	\label{fig:faandfb}
\end{figure}

Consider the following two  models $\mu^a(\theta)$ and $\mu^b(\theta)$ on the unit interval $[0,1]$, parameterized by $\theta \in (0,1)$, and with densities $f^{a}$ and $f^{b}$ given by
\[
	f^a(x,\theta)= (1-\theta) \bigl( 1+ \tfrac{1}{\theta} \bigr) (1-x)^\frac{1}{\theta} + \theta \bigl( 1 + \tfrac{1}{1-\theta} \bigr) x^\frac{1}{1-\theta},
\]
\[
	f^b(x,\theta)= \begin{cases}
	f^a(x,\theta)
	 \frac{1}{Z}\big(\one_{\{x \not\in (x_1-\frac{\delta_c}{2},x_1+\frac{\delta_c}{2})\}}+ 10^{-9} \one_{\{x \in (x_1-\frac{\delta_c}{2},x_1+\frac{\delta_c}{2})\}}\big),
	&\text{ if }\theta< 0.999,\\
	f^a(x,\theta), & \text{ if }\theta\geq 0.999,
	\end{cases}
\]
where $Z$ is a normalization constant (close to one) chosen so that $\int_{[0,1]} f^b(x,\theta)\,dx=1$.
See Figure \ref{fig:faandfb} for an illustration of these densities.

Observe that the density of model \emph{b} is that of model \emph{a} besides the small gap of width $\delta_c > 0$ created around the data point for model \emph{b} (if $\theta< 0.999$, see Figure \ref{fig:faandfb});
since the data point is fixed at $x_{1} = \tfrac{1}{2}$, the total variation distance $\rho_{\textup{TV}}\big(\mu^a(\theta),\mu^b(\theta)\big)$
between the two models is, uniformly over $\theta \in (0,1)$, bounded by a constant times $\delta_c$.
Assuming that the prior distribution on $\theta$ is the uniform distribution on $(0,1)$, observe that the prior value of the quantity of interest $\E_{\mu}[X]$ under both models (\emph{a} and \emph{b}) is approximately
$\frac{1}{2}$.
Now, when $\theta$ is close to one (zero) then the density of model \emph{a} puts most of its mass towards one (zero). Observe also that the density of model \emph{b} behaves in a similar way, with the important exception that the probability of observing the data under model \emph{b} is infinitesimally small for $\theta<0.999$.  Therefore,  for $\delta <\delta_c$, the posterior value of  the quantity of interest $\E_{\mu}[X]$ under model \emph{a} is $\frac{1}{2}$ whereas it is close to one under model \emph{b}.
Observe also that a  perturbed model \emph{c}  analogous to  \emph{b} can be constructed to lead to a posterior value close to zero.

The mechanism described here is generic and $\mu^b(\theta)$ is a simple example of what worst priors can look like after a classical Bayesian sensitivity analysis
over a class of priors specified via constraints on the TV or Prokhorov distance or the distribution of a finite number of moments.

Can these worst priors be dismissed because they depend on the data?
The problem with this argument is that, in the context of Bayesian sensitivity analysis, worst priors always depend on (or are pre-adapted to) the data.
Therefore the same argument would lead to  a dismissal of Bayesian sensitivity analysis and therefore the framework of robust Bayesian inference.
In some sense, the brittleness results reported here can  be seen as extreme occurrences of the dilation property \cite{WassermanSeidenfeld:1994} which, in robust Bayesian inference, refers to the enlargement of optimal bounds caused by the data dependence of worst priors.
Indeed, even if perturbations are quantified in  Kulback--Leibler divergence then, the local sensitivity analysis (in the sense of Fr{\'e}chet
derivatives) of posterior values \cite{GustafsonWasserman:1995}  shows infinite sensitivity as the number of data points goes to infinity (and this result is valid for the broader class of divergences that includes the Hellinger distance).

Can these worst priors be dismissed because they can ``look unrealistic'' and make the probability of observing the data very small? The problem with this  argument is that these worst priors are not ``isolated pathologies'' but directions of instability (of Bayesian conditioning) increasing with the number of data points and the complexity of the system under investigation.
We will illustrate this point with another simple example that will also show that these instabilities are the price to pay for the learning potential of Bayesian inference.

\subsection*{Learning and Robustness are Antagonistic Properties}
\label{subsec:learningvsrobustness}

In this  example we are interested in estimating $\Phi(\mu^\dagger)=\mu^\dagger[a,1]$ for some $a\in (0,1)$, where $\mu^\dagger$ is an unknown distribution on the unit interval ($\mathcal{X}=[0,1]$)  based on the observation of $n$  data point $d_1, \ldots, d_n$ up to resolution $\delta$ (i.e.\
we observe $d_i\in B_\delta(x_i)$ with $x_i\in [0,1]$ for $i=1,\ldots,n$).
Our purpose is to examine the sensitivity of the Bayesian answer to this problem with respect to the choice of a particular prior.
Consider the  model class $\mathcal{A}:=\mathcal{M}([0,1])$, and the class of priors
\[
	\Pi := \left\{ \pi \in \mathcal{M}(\mathcal{A}) \smid \E_{\mu \sim \pi} \bigl[ \E_{\mu}[X] \bigr] = m \right\}.
\]
Observe that $\Pi$ corresponds to the assumption that $\mu^\dagger$ is the realization of a random measure on $[0,1]$  whose mean, is on average $m$

As in the previous example, the finite co-dimensional class of priors $\Pi$ leads to brittleness in the sense that the least upper bound on prior values is $\mathcal{U}(\Pi)=\frac{m}{a}$,
whereas (for $\delta \ll 1/n$), the least upper bound on posterior values is the deterministic supremum of the quantity of interest (over $\mathcal{A}$), i.e.\ $\mathcal{U}(\Pi|B^n_\delta)=1$.
Furthermore, worst priors are obtained by selecting priors for which the probability of observing the data $\mu^n[B^n_\delta]$ is arbitrarily  close to zero except when $\Phi(\mu)$ is close to its deterministic supremum.

Can this brittleness be avoided by adding a uniform constraint on the probability of observing the data in the model class?
To investigate this question, let us introduce $\alpha \geq 1$ and a probability measure $\mu_0$ on $[0,1]$ with strictly positive Lebesgue density
(with $\mu_0$ being the uniform measure on $[0,1]$ as a prototypical example) and consider the (new) model class
\begin{equation}
	\label{eqadconstdata}
	\mathcal{A}(\alpha) := \left\{\mu \in \mathcal{M}[0,1] \smid \frac{1}{\alpha}\mu_0^n[B^n_\delta] \leq \mu^n[B^n_\delta] \leq \alpha\mu_0^n[B^n_\delta] \right\}
\end{equation}
and the (new) class of priors
\[
	\Pi(\alpha) := \left\{ \pi \in \mathcal{M}(\mathcal{A}(\alpha)) \smid \E_{\mu \sim \pi} \bigl[ \E_{\mu}[X] \bigr] = m \right\}
\]
where, in \eqref{eqadconstdata},  $B^{n}_{\delta}:=\prod_{i=1}^{n}{B_{\delta}(x_{i})}$ and $(x_{1},\dots,x_{n}) \in [0,1]^n$ is fixed.

Note that for the model class $\mathcal{A}(\alpha)$, the probability of observing the data is uniformly bounded from below by $\frac{1}{\alpha}\mu_0^n[B^n_\delta]$ and from above by $\alpha\mu_0^n[B^n_\delta]$.
Therefore, for $\alpha=1$,  the probability of observing the data is uniform in the model class, prior values are equal to posterior values, the method is robust but learning is impossible.
On the other hand, if $\alpha$ slightly deviates from $1$, then the calculus developed in \cite{OSS:2013} (theorems 4.8 and 4.13) gives
\begin{equation}
	\label{eqhieuhdee}
	\lim_{\delta \rightarrow 0}\mathcal{U}\big(\Pi(\alpha)|B^n_\delta\big)=\frac{1}{1+\frac{1}{\alpha^2} \frac{a-m}{m}}=\frac{m}{\frac{a}{\alpha^2}+m (1-\frac{1}{\alpha^2})}.
\end{equation}

Note that the right hand side of \eqref{eqhieuhdee} is equal to $m/a$ for $\alpha=1$ (when the probability of the data is constant on the model class)
and \emph{quickly} converges towards $1$ as $\alpha$ increases. As a numerical application observe that for $a=\frac{3}{4}$ and $m=\frac{a}{2}=\frac{3}{8}$, we have $ \lim_{\delta \rightarrow 0} \mathcal{U}\big(\Pi(\alpha)\big)=\frac{1}{2} $ and
\[
	\lim_{\delta \rightarrow 0}\mathcal{U}\big(\Pi(\alpha)|B^n_\delta\big)=\frac{1}{1+\frac{1}{\alpha^2}} .
\]
Therefore, for $\alpha=2$, we have (irrespective of the number of data points)
\[
	\lim_{\delta \rightarrow 0}\mathcal{U}\big(\Pi(2)|B^n_\delta\big)=0.8,
\]
and for $\alpha=10$, we have (irrespective of the number of data points)
\[
	\lim_{\delta \rightarrow 0}\mathcal{U}\big(\Pi(10)|B^n_\delta\big)\approx 0.99.
\]

Moreover, if $\alpha$ is derived by assuming the probability of each data point to be known up to some tolerance $\gamma$, i.e.\
if the model class $\mathcal{A}(\alpha)$ is replaced by
\begin{equation}
	\label{eqadconstdatagt}
	\mathcal{A}_\gamma := \left\{\mu \in \mathcal{M}[0,1] \smid  \frac{1}{\gamma}\mu_0[B_\delta(x_i)] \leq \mu[B_\delta(x_i)] \leq \gamma \mu_0[B_\delta(x_i)],\, i=1,\ldots,n \right\}
\end{equation}
for some $\gamma>1$, then it can be shown that
\[
	\lim_{\delta \rightarrow 0} \mathcal{U}(\Pi|B^n_\delta)=\frac{1}{1+\frac{1}{\gamma^{2n}}},
\]
which exponentially converges towards $1$ as the number $n$ of data points goes to infinity.

In conclusion,  the effects of a uniform constraint on the probability of the
data under finite information in the model class shows that learning ability comes at the price of loss in stability in the following sense: when $\alpha=1$, the data is equiprobable under all measures in the model class, posterior values are equal to prior values, the method is robust but learning is not possible. As $\alpha$ deviates from one, the learning ability increases as  robustness decreases, and when $\alpha$ is large, learning is possible but the method is brittle.

\subsection*{Qualitative Robustness and Consistency}
Since the data dependence of worst priors is inherent to classical Bayesian sensitivity analysis, one may ask whether robustness could be established under finite information by leaving the strict framework of robust Bayesian inference and computing the sensitivity of posterior conclusions independently of the specific value of the data.
Indeed, in the current classical Bayesian sensitivity analysis framework, given a class of priors $\Pi$ and the observation $d\in B^n_\delta(x)$, we compute
\[
	\sup_{\pi,\pi'\in \Pi} \Big|\E_{\mu \sim \pi}\big[\Phi(\mu)|d\in B^n_\delta(x)\big]-\E_{\mu \sim \pi'}\big[\Phi(\mu)|d\in B^n_\delta(x)\big]\Big| ,
\]
which corresponds to the \emph{sensitivity of posterior values} (given the value of the data) with respect to the particular choice of prior $\pi \in \Pi$.
Therefore, the interpretation of the brittleness mechanisms discussed above should be limited to the significance of such optimal bounds, which are not the sole measure of robustness of a Bayesian estimation.
An alternative analysis  could be to quantify the \emph{sensitivity of the distribution of posterior values}. For instance, given a class of priors $\Pi \subset \mathcal{M}(\mathcal{X})$ over a model class $\mathcal{A}\subseteq\mathcal{M}(\mathcal{X})$, the value of
\[
	\sup_{\pi,\pi'\in \Pi, \nu \in \mathcal{A}} \P_{x \sim \nu^n} \Big[ \big|\E_{\mu \sim \pi}\big[\Phi(\mu)|d\in B^n_\delta(x)\big]-\E_{\mu \sim \pi'}\big[\Phi(\mu)|d\in B^n_\delta(x)\big]\big| \geq \epsilon \Big]  .
\]
is the least upper bound on the probability that, posterior values derived from $\pi, \pi' \in \Pi$ and randomized through an admissible candidate $\nu \in \mathcal{A}$ for the distribution of the data, deviate by at least $\epsilon>0$.
This form of analysis is directly related to Hampel \cite{Hampel} and Cuevas' \cite{Cuevas} notion of \emph{qualitative robustness}, which requires  closeness in \emph{distributions of the posterior distribution}
rather than in \emph{posterior distributions}. More precisely, given a metric $\rho_2$ on
$\mathcal{M}(\mathcal{M}(\mathcal{A}))$, a \emph{qualitative} sensitivity analysis would seek to bound $\rho_2(\pi_\ast \nu^n , \pi'_\ast \nu^n )$ (over $\pi,\pi'\in \Pi$ and  $\nu \in \mathcal{A}$) where    $\pi_\ast \nu^n \in \mathcal{M}(\mathcal{M}(\mathcal{A}))$ is the distribution of the posterior distribution of the prior $\pi\in \mathcal{M}(\mathcal{A})$ when the data $d=(d_1,\ldots,d_n)$ is randomized through $\nu^n$.
If, unlike Hampel and Cuevas who require ``closeness for all $n$'', we  follow Huber \cite{HuberRonchetti:2009} and Mizera \cite{Mizera} in only requiring closeness ``for large enough $n$'' (i.e.~in the limit as the number of data points tends to infinity) then we obtain \cite{OwhadiScovel:2014} a notion of \emph{qualitative robustness} where the notion of \emph{consistency} (i.e.~the property that posterior distributions convergence towards the data generating distribution)  plays an important role.
Although consistency is primarily a frequentist notion, according to
 Blackwell and Dubins \cite{Blackwell_merging} and Diaconis and Freedman \cite{DiaconisFreedman:1986},
  consistency is equivalent to \emph{intersubjective agreement}
which means that two Bayesians will ultimately have very close predictive distributions.
Fortunately, not only are there  mild conditions which  guarantee consistency, but the posterior distributions can be shown to contract/concentrate at an exponential rate around the data generating distribution (see  \cite{Vaart2008} for rates of contraction of posterior distributions based on Gaussian process priors) and the  Bernstein--von Mises theorem goes further in providing mild conditions under which the posterior is asymptotically normal  \cite{CastilloNickl:2013, Castillo2014}. The most famous of these are Doob \cite{Doob:1949}, Le Cam and Schwartz \cite{Lecam_necessary},
 and Schwartz \cite[Thm.~6.1]{Schwartz:1965}.

Unfortunately, the conditions ensuring consistency (e.g.~the condition that the prior has Kullback--Leibler support at the parameter value generating the data\footnote{ $\pi \in \mathcal{M}(\mathcal{M}(\mathcal{X}))$ is said to have Kullback--Leibler support at $\nu\in \mathcal{M}(\mathcal{X})$ if
 $\pi\{\mu \in \mathcal{M}(\mathcal{X})\mid \int_{\mathcal{X}} \frac{d\nu}{d\mu}d\nu \leq \epsilon \}$ is strictly positive for all $\epsilon>0$}) are such that arbitrarily small (TV or Prokhorov) local perturbations of the prior distribution (near the data generating distribution) may result in consistency or non-consistency, and therefore may have large impacts on the asymptotic behavior of posterior distributions \cite{OwhadiScovel:2014}.
 A simple illustration of this mechanism is as follows \cite{OwhadiScovel:2014}. Suppose that the data generating distribution $\nu$ is at distance $\tau>0$ from the support of the prior $\pi$. Let $\pi_1$ be a prior distribution with all of its mass on/around $\nu$ (having KL support at $\nu$). Take $\pi':= (1-\epsilon)\pi+\epsilon \pi_1$. The total variation  distance from $\pi'$ to $\pi$ is bounded by $\epsilon$ which can be chosen to be arbitrarily small. Furthermore, $\pi'$ inherits the KL support of $\pi_1$ at $\nu$ and by Schwartz' consistency theorem \cite{Schwartz:1965} its posterior distribution converges (almost surely) towards a Dirac concentrated at $\nu$ as $n\rightarrow \infty$. On the other hand the distance between the support of the posterior distribution of $\pi$ and $\nu$ remains bounded by $\tau$.
 This simple example exposes a serious challenge to proving robustness in the TV metric or any weaker metric, such as those used in the convergence of MCMC.

Of course, in a parametric setting, if the parameter space $\Theta$ is compact and the model well-specified (the data is generated from a parameter in that space), then, choosing a prior satisfying Cromwell's rule (putting mass in the KL neighborhood of all parameters) ensures qualitative robustness (and the degree of robustness is  a function of how much mass is placed in each neighborhood). However, if $\Theta$ is compact and the model is misspecified then, even if the prior is nice and smooth, the mechanism discussed above suggests that it is not qualitatively robust (with a degree of non-robustness corresponding to the degree of misspecification and the prior does not need to look ``unrealistic'' to be non-qualitatively-robust).
Note also that if $\Theta$ is non-compact then the prior cannot be qualitatively robust (because no matter how small $\epsilon$ is, one can always find a neighborhood of the parameter space having mass smaller than $\epsilon$).

In a non parametric setting, consistent priors (such as the ones analysed in \cite{Vaart2008} with bounds on convergence rates) remain good/natural choices when their posterior distributions can be computed.
However, consistency and robustness are to some degree conflicting requirements \cite{Cuevas, OwhadiScovel:2014} from the point of view a numerical analyst. Consider, for instance, the problem of using a sophisticated numerical Bayesian model to predict the climate where Bayes rule is applied iteratively and posterior values become prior values for the next iteration. How do we make sure that our predictions are robust, not only with respect to the choice of prior but also with respect to numerical instabilities arising in the iterative application of the Bayes rule? The non-robustness mechanisms discussed here suggest that, unless the prior is chosen  carefully and, unless we have a  tight control on numerical instabilities/errors/approximations at each step of the iteration, our final predictions might be unstable.
Note that, oftentimes, these posterior distributions (which are later on used as prior distributions) are only approximated (e.g.\ via MCMC methods), how do we go about ensuring the stability of our method in such situations? The brittleness results discussed here suggest that having strong convergence of our MCMC method in TV would not be enough to ensure stability. Note in particular that although quantifying perturbations in KL ensures qualitative robustness it would also require controlling the convergence of the MCMC method in KL or in a stronger metric.

\subsection*{Conclusion and Perspectives}

It is possible that an analogy can be made between the brittleness and robustness properties of Bayesian inference and the numerical analysis of PDEs, for which many pathologies and also many necessary and/or sufficient stability conditions are known.
However, in contrast to conditions such as the well-known Courant--Friedrichs--Lewy condition for PDEs,  the question of the existence and of the nature of a stability condition for using Bayesian Inference under finite information remains to be answered.
Although numerical schemes that do not satisfy the CFL condition may look grossly inadequate, the existence of such perverse examples certainly does not imply the dismissal of the necessity of a stability condition. Similarly, although one may, as in the example provided in Figure \ref{fig:faandfb}, exhibit grossly perverse worst priors, the existence of such priors does not invalidate the need for a study of stability conditions for using Bayesian Inference under finite information.
The example provided in \eqref{eqhieuhdee} suggests that, in the framework of Bayesian sensitivity analysis, such a stability condition would depend on (i) how well the probability of the data is  known or  constrained in the model class (ii) the resolution at which the quantity of interest is conditioned upon the data. Note that the independence of the brittleness threshold $\delta_c$ from the number of data points $n$ in Theorem \ref{thm:1} suggests that taking $\delta$ fixed and $n\rightarrow \infty$ does not prevent brittleness in the classical Bayesian sensitivity analysis framework (it only leads to more directions of instabilities). On the other hand, for a fixed $\delta$, \eqref{eqkjhkdhjdh} suggests that brittleness results  do not persist in that same framework when the number of moment constraints $k$ (on the class of priors) is large enough. Furthermore, taking $\delta>0$ fixed (or discretizing space at a resolution $\delta>0$) enables the construction of of classes of qualitatively robust priors (to TV perturbations) that are nearly consistent as $n\rightarrow \infty$ (some degree of consistency is lost due to the discretization).
At a higher level, the mechanisms discussed here appear to suggest that robust inference (in a continuous world under finite-information) should perhaps be done with reduced/coarse models rather than highly sophisticated/complex models (with a level of ``coarseness/reduction'' depending on the available ``finite-information'').

A close inspection of some of the cases where Bayesian inference has been successful suggests the existence of a non-Bayesian feedback loop on the evaluation of its performance \cite{MayoSpanos:2004, Senn:2011, Mayo:2012}.
Therefore, one natural question is whether the missing stability condition could also be derived by exiting the strictly Bayesian framework, as proposed in \cite{Efron:2013}. One example of such an approach could be that of posterior predictive checking
 \cite{Gelman1996}, \cite[p.~159]{Gelman2004} whose rationale is to detect model mismatch by
 generating replicate data from the model, and comparing this replicate data to the original data using statistics related to the quantity of interest.

  It is natural to expect that  robustness and stability questions will increase in importance
as Bayesian methods increase in popularity due to  the availability of computational methodologies and environments to compute the posteriors. Another strong motivation for considering Bayesian methods and investigating such questions is  the complete class theorem, which, in the adversarial game theoretic setting of Decision Theory  \cite{Wald:1950} asserts that optimal statistical estimators  (leading to optimal decisions as defined by a convex loss function on a compact parameter space) live in the Bayesian class of estimators \cite{Wald:1950, Ferguson1967}.

\subsection*{Acknowledgments}
The authors gratefully acknowledge support for this work from the Air Force Office of Scientific Research under Award Number
FA9550-12-1-0389 (Scientific Computation of Optimal Statistical Estimators). The authors also thank two anonymous referees for valuable comments and suggestions.

\bibliographystyle{plain}
\bibliography{./refs}

\end{document}